\documentclass{amsart}
\usepackage{amsfonts}
\usepackage{amsmath}
\usepackage{amssymb}
\usepackage{amscd}
\usepackage{amstext}
\usepackage{amsthm}
\usepackage{enumerate}
\usepackage{color}

\theoremstyle{plain}
\newtheorem{lem}{Lemma}[section]
\newtheorem{prop}{Proposition}[section]
\newtheorem{thm}{Theorem}[section]

\newtheorem{mainthm}{Theorem}

\theoremstyle{definition}
\newtheorem{defn}{Definition}[section]

\theoremstyle{remark}
\newtheorem{ex}{Example}[section]
\newtheorem{rem}{Remark}[section]

\begin{document}
\title[Convexity properties for generalized moment maps I]
{Convexity properties \\for \\generalized moment maps I}
\author{Yasufumi Nitta}
\address{Department of Mathematics, Graduate school of science, Osaka University. 1-1 Machikaneyama, Toyonaka, Osaka 560-0043 Japan.}
\date{}
\maketitle
\thispagestyle{empty}

\begin{abstract}
We study generalized moment maps for a Hamiltonian action 
on a connected compact $H$-twisted generalized complex manifold 
introduced by Lin and Tolman and prove the convexity and 
connectedness properties of the 
generalized moment maps for a Hamiltonian torus action. 
\end{abstract}

\section{Introduction}
A notion of generalized complex structures was introduced by 
Hitchin \cite{Hi} and developed by Gualtieri \cite{Gua}. 
It provides a unifying framework for both complex and symplectic 
geometry, and a useful geometric language for understanding 
some recent development in string theory. Generalized 
K${\rm \ddot a}$hler geometry, the generalized complex geometry 
analogue of K${\rm \ddot a}$hler geometry, was introduced by 
Gualtieri, who also shows that it is essentially equivalent to that 
of a bihermitian structure, which was first discovered by physicists 
studying super-symmetric nonlinear sigma model in \cite{GHR}. 

For a group action on manifolds, notions of Hamiltonian actions 
and moment maps play a very impotant role in many geometry. It is 
an interesting and important question if there exists natural 
notions of Hamiltonian actions and moment maps. In \cite{Lin}, 
Lin and Tolman introduced notions of a Hamiltonian action 
and a generalized moment map for generalized complex geometry. 
They showed in \cite{Lin} a reduction theorem for Hamiltonian 
actions of compact Lie groups on an $H$-twisted generalized complex and K${\rm \ddot a}$hler manifold. As an application, they constructed 
explicit examples of bihermitian structures on 
$\mathbb{CP}^{n}$, Hirzebruch surfaces, the blow up of 
$\mathbb{CP}^{2}$ at arbitrarily many points, and 
other toric varieties, as well as complex Grassmannians. Their 
construction is a powerful tool for producing bihermitian structures on 
manifolds which can be produced as a symplectic reduction of 
$\mathbb{C}^{N}$. Moreover, it was shown by Kapustin and Tomasilleo 
in \cite{KT} that the mathematical notion of Hamiltonian actions 
on a generalized K${\rm \ddot a}$hler manifold corresponds exactly 
to the physiccal notion of general $(2, 2)$ gauged sigma models 
with $3$-form fluxes. 

Convexity and connectedness properties for moment maps of 
Hamiltonian torus actions on a connected compact symplectic manifold was 
shown by Atiyah \cite{At} and Guillemin and Sternberg \cite{GS}. 
In the present paper, we study Hamiltonian torus actions on a 
connected compact $H$-twisted generalized complex manifold 
and prove the convexity and 
connectedness properties of a generalized moment map for 
Hamiltonian torus actions. The main result is stated below. 
The detailed notations and definitions are in section 2 
and section 3. 
\begin{mainthm}
Let an $m$-dimensional torus $T^{m}$ act on a connected compact 
$H$-twisted generalized complex manifold $(M, \mathcal{J})$ in a 
Hamiltonian way with a generalized moment map 
$\mu : M \longrightarrow \mathfrak{t}^{*}$ and a moment one 
form $\alpha \in \Omega^{1}(M; \mathfrak{t}^{*})$. 
Then: 
\begin{enumerate}
\item the levels of $\mu$ are connected; 
\item the image of $\mu$ is convex; 
\item the fixed points of the action form a finite 
union of connected submanifolds $C_{1}, \cdots, C_{N}$:
\begin{equation*}
{\rm Fix}(T^{m}) = \bigcup_{i=1}^{N}C_{i}. 
\end{equation*}
On each component the generalized moment map 
$\mu$ is constant: $\mu(C_{i}) = \{a_{i}\}$, and the image of 
$\mu$ is the convex hull of the images $a_{1}, \cdots, a_{N}$ 
of the fixed points of the action, that is, 
\begin{equation*}
\mu(M) = 
\left\{ \sum_{i=1}^{N}\lambda_{i}a_{i}\ |\ \sum_{i=1}^{N}\lambda_{i}=1, \lambda_{i} \geq 0 \right\}. 
\end{equation*}
\end{enumerate}
\end{mainthm}
This paper is organized as follows. In section 2 we briefly review 
of the theory of generalized complex structures and generalized 
K${\rm \ddot a}$hler structures. In section 3 we introduce the 
notion of generalized moment maps for Hamiltonian actions on a 
generalized complex manifold and prove that the generalized 
moment map has a property of a Bott-Morse function. At the last section, 
we shall give a proof of Theorem A. 
\section{Generalized complex structures}
First we recall the basic theory of generalized complex 
structures; see \cite{Gua} for the details. 

Given a closed $3$-form $H$ on an $n$-dimensional manifold $M$, 
we define the $H$-twisted Courant bracket of sections of 
the sum $T \oplus T^{*}$ of the tangent and cotangent bundles by 
\begin{equation*}
[X + \xi, Y + \eta]_{H} = 
[X, Y] + \mathcal{L}_{X}\eta - \mathcal{L}_{Y}\xi 
- \frac{1}{2}d\left( \eta(X) - \xi(Y) \right) + i_{Y}i_{X}H, 
\end{equation*}
where $\mathcal{L}_{X}$ denotes the Lie derivative along a vector 
field $X$. The vector bundle $T \oplus T^{*}$ is also endowed 
with a natural inner product of signature $(n, n)$: 
\begin{equation*}
\langle X + \xi, Y + \eta \rangle = \frac{1}{2}(\eta(X) + \xi(Y)). 
\end{equation*}

\begin{defn}
Let $M$ be a manifold, and $H$ be a closed $3$-form on $M$. 
A generalized almost complex structure on $M$ is a complex 
structure $\mathcal{J}$ on the bundle $T \oplus T^{*}$ 
which preserves the natural inner product. 
If sections of the $\sqrt{-1}$-eigenspace of $\mathcal{J}$ is closed 
under the $H$-twisted Courant bracket, $\mathcal{J}$ is called an 
$H$-twisted generalized complex structure. If $H = 0$, we call it 
simply a generalized complex structure. 
\end{defn}
An $H$-twisted generalized complex structure can be fully 
described in terms of its $\sqrt{-1}$-eigenspace $L$, 
which is a maximal isotropic subspace of 
$(T \oplus T^{*}) \otimes \mathbb{C}$ satisfying 
$L \cap \bar L = \{ 0 \}$ and to be closed under 
the $H$-twisted Courant bracket. 

Let $\pi : (T \oplus T^{*}) \otimes \mathbb{C} \longrightarrow 
T \otimes \mathbb{C}$ be the natural projection. The type of an 
$H$-twisted generalized complex structure $\mathcal{J}$ is 
the codimension of $\pi(L)$ in $T \otimes \mathbb{C}$, where $L$ is 
the $\sqrt{-1}$-eigenspace of $\mathcal{J}$. 

\begin{ex}[Complex structures (type $n$)]
Let $J$ be a usual complex structure on a $2n$-dimensional manifold $M$. 
Consider the endomorphism 
\begin{equation*}
\mathcal{J}_{J}=
\left(
\begin{array}{ccc}
 J  & 0      \\
 0  & -J^{*} \\
\end{array}
\right), 
\end{equation*}
where the matrix is written with respect to the direct sum 
$T \oplus T^{*}$. Then $\mathcal{J}_{J}$ is a 
generalized complex structure of type $n$ on $M$; 
the $\sqrt{-1}$-eigenspace of $\mathcal{J}_{J}$ is 
$L_{J} = T_{1,0} \oplus T^{0,1}$, where $T_{1,0}$ is the 
$\sqrt{-1}$-eigenspace of $J$. 
\end{ex}

\begin{ex}[Symplectic structures (type 0)]
Let $\omega$ be a symplectic structure on a $2n$-dimensional manifold 
$M$, viewed as a skew-symmetric isomorphism 
$\omega :T \longrightarrow T^{*}$ via the interior product 
$X \mapsto i_{X}\omega$. Consider the endomorphism 
\begin{equation*}
\mathcal{J}_{\omega}=
\left(
\begin{array}{ccc}
 0       & -\omega^{-1} \\
 \omega  & 0            \\
\end{array}
\right). 
\end{equation*}
Then $\mathcal{J}_{\omega}$ is a 
generalized complex structure of type $0$ on $M$; 
the $\sqrt{-1}$-eigenspace of $\mathcal{J}_{\omega}$ is 
\begin{equation*}
L_{\omega} = \{ X - \sqrt{-1}i_{X}\omega |\ X \in T \otimes \mathbb{C} \}.
\end{equation*}
\end{ex}

\begin{ex}[$B$-field shift]
Let $(M, \mathcal{J})$ be an $H$-twisted generalized complex manifold 
and $B \in \Omega^{2}(M)$ be a closed $2$-form on $M$. Then the endomorphism 
\begin{equation*}
\mathcal{J}_{B}=
\left(
\begin{array}{ccc}
 1 & 0 \\
 B & 1 \\
\end{array}
\right)
\mathcal{J}
\left(
\begin{array}{ccc}
 1   & 0 \\
 -B  & 1 \\
\end{array}
\right)
\end{equation*}
is also an $H$-twisted generalized complex structure. It is called the 
$B$-field shift of $\mathcal{J}$. The type of $\mathcal{J}_{B}$ coincides 
with that of $\mathcal{J}$. The $\sqrt{-1}$ eigenspace $L_{B}$ of 
$\mathcal{J}_{B}$ can be written by 
\begin{equation*}
L_{B} = \{ X + f + i_{X}B\ |\ X + f \in L \}, 
\end{equation*}
where $L$ is the $\sqrt{-1}$ eigenspace of $\mathcal{J}$. 
\end{ex}
The type of an $H$-twisted generalized complex structure is not required to 
be constant along the manifold, and it may jump along loci. 
Gualtieri constructed a generalized complex structure on 
$\mathbb{CP}^{2}$ which is type 2 along a cubic curve and type 0 
outside the cubic curve. The detailed construction can be seen in 
\cite{Gua}. 

Next we briefly review the notion of 
$H$-twisted generalized K${\rm \ddot a}$hler structures. 

\begin{defn}
Let $M$ be a manifold, and $H$ be a closed $3$-form on $M$. 
an $H$-twisted generalized K${\rm \ddot a}$hler structure on $M$ 
is a pair of commuting $H$-twisted generalized complex structures 
$(\mathcal{J}_{1}, \mathcal{J}_{2})$ so that 
$\mathcal{G} = -\mathcal{J}_{1} \mathcal{J}_{2}$ is 
a positive definite metric, 
that is, $\mathcal{G}^{2} = {\rm id}$, $\mathcal{G}$ preserves the natural inner product, and 
$\mathcal{G}(X + \xi, Y + \eta) := 
\langle \mathcal{G}(X + \xi), X + \xi \rangle > 0$ 
for all non-zero $X + \xi \in T \oplus T^{*}$. 
\end{defn}

\begin{ex}\
\begin{enumerate}
\item
Let $(M, g, J)$ be a K${\rm \ddot a}$hler manifold and $\omega = gJ$ 
be the K${\rm \ddot a}$hler form. By examples above, $J$ and $\omega$ 
induce generalized complex structures $\mathcal{J}_{J}$ 
and $\mathcal{J}_{\omega}$, respectively. Moreover, $\mathcal{J}_{J}$ 
commutes with $\mathcal{J}_{\omega}$ and 
\begin{equation*}
\mathcal{G} = -\mathcal{J}_{J}\mathcal{J}_{\omega} = 
\left(
\begin{array}{ccc}
 0 & g^{-1} \\
 g & 0 \\
\end{array}
\right)
\end{equation*}
is a positive definite metric on $T \oplus T^{*}$. Hence 
$(\mathcal{J}_{J}, \mathcal{J}_{\omega})$ is a generalized 
K${\rm \ddot a}$hler structure on $M$. 
\item 
Let $(\mathcal{J}_{1}, \mathcal{J}_{2})$ be an $H$-twisted generalized 
K${\rm \ddot a}$hler structure, and $B$ be a closed $2$-form on $M$. 
Then $((\mathcal{J}_{1})_{B}, (\mathcal{J}_{2})_{B})$ is also an 
$H$-twisted generalized K${\rm \ddot a}$hler structure. It is called 
the $B$-field shift of $(\mathcal{J}_{1}, \mathcal{J}_{2})$. 
\end{enumerate}
\end{ex}

In \cite{Gua}, a characterization of $H$-twisted generalized 
K${\rm \ddot a}$hler pairs was given in terms of Hermitian 
geometry, which is represented below. 

\begin{thm}[M. Gualtieri, \cite{Gua}]
For each $H$-twisted generalized K${\rm \ddot a}$hler structure 
$(\mathcal{J}_{1}, \mathcal{J}_{2})$, there exists unique $2$-form $b$, 
Riemannian metric $g$, and two orthogonal complex structures $J_{\pm}$ 
such that 
\begin{equation*}
\mathcal{J}_{1,2} = \frac{1}{2}
\left(
\begin{array}{ccc}
 1 & 0 \\
 b & 1 \\
\end{array}
\right)\left(
\begin{array}{ccc}
 J_{+} \pm J_{-}           & -(\omega^{-1}_{+} \mp \omega^{-1}_{-}) \\
 \omega_{+} \mp \omega_{-} & -(J^{*}_{+} \pm J^{*}_{-})            \\
\end{array}
\right)\left(
\begin{array}{ccc}
 1  & 0 \\
 -b & 1 \\
\end{array}
\right), 
\end{equation*}
where $\omega_{\pm} = gJ_{\pm}$ satisfy 
\begin{equation}\label{torsion}
d^{c}_{-}\omega_{-} = -d^{c}_{+}\omega_{+} = H + db. 
\end{equation}
Conversely, any quadruple $(g, b, J_{\pm})$ satisfying condition 
$(\ref{torsion})$ defines an $H$-twisted generalized 
K${\rm \ddot a}$hler structure. 
\end{thm}

Every $H$-twisted generalized complex manifold may not admit an 
$H$-twisted generalized K${\rm \ddot a}$hler structure. However, 
following lemma claims that every $H$-twisted generalized complex 
manifold admits a generalized almost K${\rm \ddot a}$hler structure. 
This is a generalized complex geometry analogue of the fact 
that a symplectic manifold admits an almost complex structure 
which is compatible with the symplectic structure. 

\begin{lem}\label{compatible1}
Let $(M, \mathcal{J})$ be an $H$-twisted generalized complex 
manifold. Then there exists a generalized almost complex structure 
$\mathcal{J}^{\prime}$ which is compatible with $\mathcal{J}$, that is, 
$\mathcal{J}^{\prime}$ is a generalized almost complex structure which 
commutes with $\mathcal{J}$ and 
$\mathcal{G} = -\mathcal{J}\mathcal{J}^{\prime}$ is a positive 
definite metric. 
\end{lem}

\begin{proof}
Choose a Riemannian metric $g$ on $M$ and put 
\begin{equation*}
\tilde{\mathcal{G}} = 
\left(
\begin{array}{ccc}
 0 & g^{-1} \\
 g & 0 \\
\end{array}
\right). 
\end{equation*}
Then $\tilde{\mathcal{G}}$ is a postitive definite metric on $T \oplus T^{*}$. 
Now we define a symplectic structure $\mathcal{W}$ on $T \oplus T^{*}$ 
by 
\begin{equation*}
\mathcal{W}(X + \xi, Y + \eta) = 
\langle \mathcal{J}(X + \xi), Y + \eta \rangle . 
\end{equation*}
Since $\tilde{\mathcal{G}}$ and $\mathcal{W}$ are non-degenerate, 
there exists an endomorphism $\mathcal{A}$ on $T \oplus T^{*}$ 
which satisfies 
\begin{equation*}
\mathcal{W}(X + \xi, Y + \eta) = 
\tilde{\mathcal{G}}(\mathcal{A}(X + \xi), Y + \eta) 
\end{equation*}
for all $X + \xi, Y + \eta \in T \oplus T^{*}$. 
The map $\mathcal{A}$ is skew-symmetric because 
\begin{eqnarray*}
\tilde{\mathcal{G}}(\mathcal{A}^{*}(X + \xi), Y + \eta) &=& 
\tilde{\mathcal{G}}(X + \xi, \mathcal{A}(Y + \eta)) = 
\tilde{\mathcal{G}}(\mathcal{A}(Y + \eta), X + \xi) \\ 
&=& \mathcal{W}(Y + \eta, X + \xi) = 
-\mathcal{W}(X + \xi, Y + \eta) \\
&=& \tilde{\mathcal{G}}(-\mathcal{A}(X + \xi), Y + \eta), 
\end{eqnarray*}
where $\mathcal{A}^{*}$ denotes the adjoint operater of $\mathcal{A}$ 
with respect to the positive difinite metric $\tilde{\mathcal{G}}$. 
Moreover since $\mathcal{A}$ is invertible, $\mathcal{A}\mathcal{A}^{*}$ 
is symmetric and positive, that is, 
$(\mathcal{A}\mathcal{A}^{*})^{*} = \mathcal{A}\mathcal{A}^{*}$ and 
\begin{equation*}
\tilde{\mathcal{G}}(\mathcal{A}\mathcal{A}^{*}(X + \xi), X + \xi ) > 0
\end{equation*}
for all non-zero $X + \xi \in T \oplus T^{*}$. 
So we can define $\sqrt{\mathcal{A}\mathcal{A}^{*}}$ the square root 
of $\mathcal{A}\mathcal{A}^{*}$. Of course 
$\sqrt{\mathcal{A}\mathcal{A}^{*}}$ is also symmetric and positive definite. 
Let $\mathcal{J}^{\prime}$ be an endomorphism on $T \oplus T^{*}$ defined by 
$\mathcal{J}^{\prime} = (\sqrt{\mathcal{A}\mathcal{A}^{*}})^{-1}\mathcal{A}$. 
Since $\mathcal{A}$ commutes with $\sqrt{\mathcal{A}\mathcal{A}^{*}}$, 
$\mathcal{J}^{\prime}$ also commutes with $\mathcal{A}$ and 
$\sqrt{\mathcal{A}\mathcal{A}^{*}}$. Hence $\mathcal{J}^{\prime}$ is 
orthogonal with respect to $\tilde{\mathcal{G}}$ and 
$(\mathcal{J}^{\prime})^{2} = -{\rm id}$. 

By the definition of $\mathcal{A}$ we have 
$\mathcal{A}\mathcal{J} = -\mathcal{J}\mathcal{A}^{-1}$, and 
$\sqrt{\mathcal{A}\mathcal{A}^{*}}\mathcal{J} = 
\mathcal{J}(\sqrt{\mathcal{A}\mathcal{A^{*}}})^{-1}$. 
Therefore we have 
\begin{eqnarray*}
\mathcal{J}\mathcal{J}^{\prime} 
&=& \mathcal{J}(\sqrt{\mathcal{A}\mathcal{A}^{*}})^{-1}\mathcal{A} \\
&=& -\sqrt{\mathcal{A}\mathcal{A}^{*}}\mathcal{A}^{-1}\mathcal{J} \\
&=& -(\mathcal{J}^{\prime})^{-1}\mathcal{J} \\
&=& \mathcal{J}^{\prime}\mathcal{J}, 
\end{eqnarray*}
in particular $\mathcal{J}^{\prime}$ commutes 
with $\mathcal{J}$. Moreover, for each 
$X + \xi, Y + \eta \in T \oplus T^{*}$, we have 
\begin{eqnarray*}
\langle \mathcal{J}^{\prime}(X + \xi), \mathcal{J}^{\prime}(Y + \eta) \rangle 
&=& -\mathcal{W}(\mathcal{J}\mathcal{J}^{\prime}(X + \xi), \mathcal{J}^{\prime}
(Y + \eta)) \\
&=& -\tilde{\mathcal{G}}(\mathcal{A}\mathcal{J}\mathcal{J}^{\prime}(X + \xi), 
\mathcal{J}^{\prime}(Y + \eta)) \\
&=& -\tilde{\mathcal{G}}(\mathcal{J}^{\prime}\mathcal{A}\mathcal{J}(X + \xi), 
\mathcal{J}^{\prime}(Y + \eta)) \\
&=& -\tilde{\mathcal{G}}(\mathcal{A}\mathcal{J}(X + \xi), Y + \eta) \\
&=& -\mathcal{W}(\mathcal{J}(X + \xi), Y + \eta) \\
&=& \langle X + \xi, Y + \eta \rangle. \\
\end{eqnarray*}
Hence $\mathcal{J}^{\prime}$ is a generalized almost 
complex structure on $M$ which commutes $\mathcal{J}$. 
Finally $\mathcal{G} := -\mathcal{J}\mathcal{J}^{\prime}$ is a 
positive definite metric on $T \oplus T^{*}$ 
since 
$\mathcal{G} = 
\tilde{\mathcal{G}}\sqrt{\mathcal{A}\mathcal{A}^{*}}$. 
This completes the proof. 
\end{proof}

If $\mathcal{J}^{\prime}$ is a generalized almost complex 
structure which is compatible with an $H$-twisted generalized 
complex structure $\mathcal{J}$, then we can apply the argument 
of Gualtieri in \cite{Gua} and construct a Riemannian metric $g$, 
a $2$-form $b$, and two orthogonal almost complex structures 
$J_{\pm}$ which satisfy the equation
\begin{equation}\label{decomp}
\mathcal{J} = \frac{1}{2}
\left(
\begin{array}{ccc}
 1 & 0 \\
 b & 1 \\
\end{array}
\right)\left(
\begin{array}{ccc}
 J_{+} + J_{-}           & -(\omega^{-1}_{+} - \omega^{-1}_{-}) \\
 \omega_{+} - \omega_{-} & -(J^{*}_{+} + J^{*}_{-})            \\
\end{array}
\right)\left(
\begin{array}{ccc}
 1  & 0 \\
 -b & 1 \\
\end{array}
\right). 
\end{equation}
Of course, $J_{+}$ and $J_{-}$ are not integrable in general. 
\section{Hamiltonian action on generalized manifolds}
In this section we introduce the definition of 
Hamiltonian actions on $H$-twisted generalized complex manifolds 
given in \cite{Lin}. 

\begin{defn}[Y. Lin and S. Tolman, \cite{Lin}]
Let a compact Lie group $G$ with its Lie algebra $\mathfrak{g}$ 
act on an $H$-twisted generalized complex manifold $(M, \mathcal{J})$ 
preserving $\mathcal{J}$, where $H \in \Omega^{3}(M)^{G}$ is a closed 
$3$-form. The action of $G$ is said to be Hamiltonian 
if there exists a $G$-equivariant 
smooth function $\mu : M \longrightarrow \mathfrak{g}^{*}$, 
called the generalized moment map, and a $\mathfrak{g}^{*}$-valued 
one form $\alpha \in \Omega^{1}(M, \mathfrak{g}^{*})$, 
called the moment one form such that 
\begin{itemize}
\item $\xi_{M} - \sqrt{-1}(d\mu^{\xi} + \sqrt{-1}\alpha^{\xi})$ lies in 
$L$ for all $\xi \in \mathfrak{g}$, where $\xi_{M}$ denotes 
the induced vector field on $M$ and $L \subset (T \oplus T^{*}) 
\otimes \mathbb{C}$ denotes the $\sqrt{-1}$-eigenspace of $\mathcal{J}$, and 
\item $i_{\xi_{M}}H = d\alpha^{\xi}$ for all $\xi \in \mathfrak{g}$. 
\end{itemize}
\end{defn}

\begin{ex}\ 
\begin{enumerate}
\item Let $G$ act on a symplectic manifold $(M, \omega)$ preserving 
$\omega$, and $\mu : M \longrightarrow \mathfrak{g}^{*}$ be an 
usual moment map, that is, $\mu$ is $G$-equivariant and 
$i_{\xi_{M}}\omega = d\mu^{\xi}$ for all $\xi \in \mathfrak{g}$. 
Then $G$ also preserves $\mathcal{J}_{\omega}$, $\mu$ is also a 
generalized moment map, and $\alpha = 0$ is a moment one form for 
this action. Hence the $G$-action on $(M, \mathcal{J}_{\omega})$ 
is Hamiltonian. 
\item Let $(M, J)$ be a complex manifold and $G$ act on 
$(M, \mathcal{J}_{J})$ in a Hamiltonian way. Then $G$ also 
preserves the original complex structure $J$. 
Since $L_ {J} = T_{1, 0} \oplus T^{0,1}$ and $\xi_{M} \in \pi(L_{J})$, 
so we have $\xi_{M} = 0$ for all $\xi \in \mathfrak{g}$. Thus if $G$ is 
connected, the $G$-action on $M$ must be trivial. 
\item Let $G$ act on an $H$-twisted generalized complex manifold 
$(M, \mathcal{J})$ with a generalized moment map $\mu$ and a moment 
one form $\alpha$. If $B \in \Omega^{2}(M)^{G}$ is closed, then $G$ 
acts on $M$ preserving the $B$-field shift of $\mathcal{J}$ with 
generalized moment map $\mu$ and moment one form $\alpha^{\prime}$, where 
$(\alpha^{\prime})^{\xi} = \alpha^{\xi} + i_{\xi_{M}}B$ for all 
$\xi \in \mathfrak{g}$. 
\end{enumerate}
\end{ex}

By the definition and examples above, we can say that the notion of 
generalized moment maps is a generalization of the notion of moment maps 
in symplectic geometry. Generalized moment maps are studied by Lin and 
Tolman in \cite{Lin}. They showed in \cite{Lin} that a reduction 
theorem for Hamiltonian actions of compact Lie groups on 
an $H$-twisted generalized complex and K${\rm \ddot a}$hler manifold 
holds. 

Now we can state Theorem A in Introduction. Before we begin a 
proof, we prove a notable property of generalized moment maps. 
At first we prove following lemmata. 

\begin{lem}\label{compatible2}
Let a compact Lie group $G$ act on an $H$-twisted generalized conplex 
manifold $(M, \mathcal{J})$ preserving $\mathcal{J}$. Then there exists a 
$G$-invariant generalized almost complex structure which is compatible 
with $\mathcal{J}$. 
\end{lem}

\begin{proof}
Choose a $G$-invariant Riemannian metric $g$ on $M$ and put 
\begin{equation*}
\mathcal{G} = 
\left(
\begin{array}{ccc}
 0 & g^{-1} \\
 g & 0 \\
\end{array}
\right). 
\end{equation*}
Then $\mathcal{G}$ is a $G$-invariant positive definite metric on 
$T \oplus T^{*}$. 
Let $\mathcal{A}$ be an endomorphism on $T \oplus T^{*}$ defined by 
$\mathcal{A} = \mathcal{G}^{-1}\mathcal{J}$. 
Since $\mathcal{G}$ and $\mathcal{J}$ are $G$-invariant, 
so $\mathcal{A}$ is also $G$-invariant. Now if we define 
\begin{equation*}
\mathcal{J}^{\prime} = (\sqrt{\mathcal{A}\mathcal{A}^{*}})^{-1}\mathcal{A}, 
\end{equation*}
then $\mathcal{J}^{\prime}$ is a generalized almost complex structure 
on $M$ which is compatible with $\mathcal{J}$. Moreover since 
$\mathcal{A}$ is $G$-invariant, $\mathcal{J}^{\prime}$ is also $G$-invariant. 
This completes the proof. 
\end{proof}

\begin{lem}\label{Fix}
Let an $m$-dimensional torus $T^{m}$ act on an $H$-twisted generalized conplex 
manifold $(M, \mathcal{J})$ in a Hamiltonian way with a generalized 
moment map $\mu$ and a moment one form $\alpha$. Then for an arbitrary 
subtorus $G \subset T^{m}$ the fixed point set of $G$-action 
\[
{\rm Fix}(G) = \{ p \in M\ |\ \theta \cdot p = p\ (\forall \theta \in G) \}
\]
is an even dimensional submanifold of $M$. 
\end{lem}

\begin{proof}
Choose a $G$-invariant generalized almost complex structure 
$\mathcal{J}^{\prime}$ which is compatible with $\mathcal{J}$. 
Then there exists a Riemannian metric $g$, a $2$-form $b$, and 
two orthogonal almost complex structures $J_{\pm}$ which satisfies 
the equation (\ref{decomp}). 
Since $\mathcal{J}$ and $\mathcal{J}^{\prime}$ are $G$-invariant, 
so $g$ and $J_{\pm}$ are also $G$-invariant. 
For each $p \in {\rm Fix}(G)$ and $\theta \in G$, 
the differential of the action of $\theta$ at $p$, 
\begin{equation*}
(\theta_{*})_{p} : T_{p}M \longrightarrow T_{p}M, 
\end{equation*}
preserves the almost complex structures $J_{\pm}$. In addition, since 
$G$-action preserves the metric $g$, the exponential mapping 
$\exp_{p} : T_{p}M \longrightarrow M$ with respect to the metric $g$ 
is equivariant, that is, 
\begin{equation*}
\exp_{p}((\theta_{*})_{p}X) = \theta \cdot \exp_{p}X 
\end{equation*}
for any $\theta \in G$ and $X \in T_{p}M$. This concludes that the 
fixed point of the action of $\theta$ near $p$ correspond to 
the fixed point of $(\theta_{*})_{p}$ on $T_{p}M$ by the exponential 
mapping, that is, 
\begin{equation*}
T_{p}{\rm Fix}(G) = \bigcap_{\theta \in G}\ker(1 - (\theta_{*})_{p}). 
\end{equation*}
Since $(\theta_{*})_{p} \circ J_{\pm} = J_{\pm} \circ (\theta_{*})_{p}$, 
the eigenspace with eigenvalue $1$ of $(\theta_{*})_{p}$ is 
invariant under $J_{\pm}$, and is therefore an almost complex 
subspace. In particular, $T_{p}{\rm Fix}(G)$ is even dimensional. 
\end{proof}

We remark in the proof of Proposition \ref{Fix} 
that the fixed point set ${\rm Fix}(G)$ is an almost 
complex submanifold with respect to $J_{\pm}$. In particular, 
we see $\omega_{\pm}$ is non-degenerate on ${\rm Fix}(G)$. Moreover, 
it is known that ${\rm Fix}(G)$ is a generalized complex 
submanifold when $H = 0$ (see \cite{Lin2} for the details). 

The following proposition claims that a generalized moment map has a 
property of a Bott-Morse function. This fact plays a crucial role 
in the proof of Theorem A. 

\begin{prop}\label{BM}
Let an $m$-dimeansional torus $T^{m}$ act on a compact 
$H$-twisted generalized conplex manifold $(M, \mathcal{J})$ 
in a Hamiltonian way with a generalized 
moment map $\mu$ and a moment one form $\alpha$. Then 
$\mu^{\xi}$ is a Bott-Morse function with even index and coindex 
for all $\xi \in \mathfrak{t}$. 
\end{prop}

\begin{proof}
Let $\xi \in \mathfrak{t}$ and $T^{\xi}$ denote the subtorus of 
$T^{m}$ which is generated by $\xi$. First we shall prove 
that the critical set 
\begin{equation*}
{\rm Crit}(\mu^{\xi}) = \{ p \in M\ |\ (d\mu^{\xi})_{p} = 0 \} 
\end{equation*}
coincides with the fixed point set of $T^{\xi}$-action ${\rm Fix}(T^{\xi})$. 
Choose a $T^{m}$-invariant generalized almost 
complex structure $\mathcal{J}^{\prime}$ which is 
compatible with $\mathcal{J}$. Then $\mathcal{J}$ can be written 
by the form of the equation (\ref{decomp}) by corresponding 
quadruple $(g, b, J_{\pm})$. Note that the metric $g$ and 
orthogonal almost complex structures $J_{\pm}$ are all $T^{m}$-invariant. 

Since $\xi_{M} - \sqrt{-1}(d\mu^{\xi} + \sqrt{-1}\alpha^{\xi}) \in L$ 
by the definition of Hamiltonian actions, so $(d\mu^{\xi})_{p} = 0$ 
implies $p \in {\rm Fix}(T^{\xi})$. 
In particular we obtain ${\rm Crit}(\mu^{\xi}) \subset {\rm Fix}(T^{\xi})$. 
On the other hand, since 
${\rm Fix}(T^{\xi}) = \{ p \in M\ |\ (\xi_{M})_{p} = 0\}$ 
so we can view $\mu^{\xi}$ locally as an imaginary part of a 
pseudoholomorphic function on an almost complex manifold 
$({\rm Fix}(T^{\xi}), J_{\pm})$. By applying the 
principle of the maximum and compactness of ${\rm Fix}(T^{m})$, 
we see that $\mu^{\xi}$ is constant on each connected component of 
${\rm Fix}(T^{\xi})$. Moreover the gradient of $\mu^{\xi}$ with 
respect to the metric $g$ is 
tangent to ${\rm Fix}(T^{\xi})$ 
on there because $g$ and $\mu^{\xi}$ are $T^{\xi}$-invariant. 
This shows that ${\rm Fix}(T^{\xi}) \subset {\rm Crit}(\mu^{\xi})$, 
and hence we obtain ${\rm Crit}(\mu^{\xi}) = {\rm Fix}(T^{\xi})$. 
In particular, ${\rm Crit}(\mu^{\xi})$ is 
an even dimensional submanifold of $M$. 

To prove $\mu^{\xi}$ is a Bott-Morse function, we shall calculate the 
Hessian of $\mu^{\xi}$ on ${\rm Crit}(\mu^{\xi})$. First note that 
the induced vector field 
$\xi_{M}$ can be written by 
\begin{equation*}
\xi_{M} = \frac{1}{2}
\Big( \omega_{+}^{-1}(d\mu^{\xi}) - \omega_{-}^{-1}(d\mu^{\xi}) \Big). 
\end{equation*}
Let $\nabla$ be the Riemannian connection 
with respect to $g$, and $\nabla^{2}\mu^{\xi}$ denotes the 
Hessian of $\mu^{\xi}$. Then by an easy calculation we have 
\begin{eqnarray*}
g(\nabla^{2}\mu^{\xi}(Y), Z) &=& 
Y(Z\mu^{\xi}) - \nabla_{Y}Z(\mu^{\xi}) \\
&=& Y(g(J_{\pm}\xi_{M}^{\pm}, Z )) - g(J_{\pm}\xi_{M}^{\pm}, \nabla_{Y}Z) \\
&=& g(\nabla_{Y}J_{\pm}\xi_{M}^{\pm}, Z) \\
&=& g((\nabla_{Y}J_{\pm})\xi_{M}^{\pm}, Z) + g(J_{\pm}(\nabla_{Y}\xi_{M}^{\pm}), Z), 
\end{eqnarray*}
where $\xi_{M}^{\pm} = \omega_{\pm}^{-1}(d\mu^{\xi}) = 
-J_{\pm}g^{-1}(d\mu^{\xi})$. Thus for each $p \in {\rm Crit}(\mu^{\xi})$ 
we have 
\begin{equation*}
(\nabla^{2}\mu^{\xi})_{p} = J_{\pm}(\nabla_{Y_{p}}\xi_{M}^{\pm}) 
\end{equation*}
because $(\xi_{M}^{\pm})_{p} = 0$. 
Let $(L_{\xi})_{p}$ be an endomorphism on $T_{p}M$ 
defined by 
\begin{equation*}
(L_{\xi})_{p}(Y) = [\xi_{M}, Y]_{p} = -\nabla_{Y_{p}}\xi_{M}. 
\end{equation*}
Then since 
$\xi_{M} = \frac{1}{2}
\left( \xi_{M}^{+} - \xi_{M}^{-} \right)$, 
$(L_{\xi})_{p}$ can be written by 
\begin{equation*}
(L_{\xi})_{p} = -\frac{1}{2}(J_{+} - J_{-})(\nabla^{2}\mu^{\xi})_{p}. 
\end{equation*}

Now we prove that 
$T_{p}{\rm Crit}(\mu^{\xi}) = \ker (\nabla^{2}\mu^{\xi})_{p}$. 
Since ${\rm Crit}(\mu^{\xi})$ is a submanifold of $M$, it is easy to see 
$T_{p}{\rm Crit}(\mu^{\xi}) \subset \ker (\nabla^{2}\mu^{\xi})_{p}$. 
So we may only show that 
$\ker (\nabla^{2}\mu^{\xi})_{p} \subset T_{p}{\rm Crit}(\mu^{\xi})$. 
At first we have 
$\ker (\nabla^{2}\mu^{\xi})_{p} \subset \ker (L_{\xi})_{p}$ 
by the calculation above. 
If we identify $(L_{\xi})_{p}$ with a vector field on $T_{p}M$, 
the one parameter family of diffeomorphism 
$\{ (\exp t\xi _{*})_{p} \}_{t \in \mathbb{R}}$ on $T_{p}M$ coincides with 
$\{ \exp t(L_{\xi})_{p} \}_{t \in \mathbb{R}}$. So $\ker (L_{\xi})_{p}$ 
coincides with the fixed point set of 
$\{ (\exp t\xi _{*})_{p} \}_{t \in \mathbb{R}}$. 
Hence we have 
\begin{equation*}
\ker (\nabla^{2}\mu^{\xi})_{p} 
\subset \ker (L_{\xi})_{p} 
= \bigcap_{\theta \in T^{\xi}}\ker (1 - (\theta_{*})_{p}) 
= T_{p}{\rm Fix}(T^{\xi})
= T_{p}{\rm Crit}(\mu^{\xi}), 
\end{equation*}
and this shows that 
$T_{p}{\rm Crit}(\mu^{\xi}) = \ker (\nabla^{2}\mu^{\xi})_{p}$. 
In particular, $\mu^{\xi}$ is a Bott-Morse function. 

Finally, by an easy calculation we see that $(\nabla^{2}\mu^{\xi})_{p}$ 
commutes with $J_{+} - J_{-}$ for all $p \in {\rm Crit}(\mu^{\xi})$. 
So we can define a non-degenerate $2$-form on each non-zero 
eigenspace of $(\nabla^{2}\mu^{\xi})_{p}$ by $g(J_{+}-J_{-})$. 
Hence each non-zero eigenspace of $(\nabla^{2}\mu^{\xi})_{p}$ is 
even dimensional, in particular the index and coindex of 
the critical manifold are even. 
\end{proof}

\begin{rem}
If $M$ is noncompact, then the generalized moment map is not a 
Bott-Morse function in general. Indeed, if we consider 
a trivial circle action on a complex manifold $(M, J)$, 
then the imaginary part of an arbitrary holomorphic 
function is a generalized moment map for this action. 
\end{rem}

\section{Proof of Theorem A}
We shall prove Theorem A in this section. 
This proof involves induction over $m = \dim T^{m}$. 
Consider the statements: 
\begin{eqnarray*}
A_{m} &:& 
\text{"the level sets of $\mu$ are connected, for any $T^{m}$-action",} \\
B_{m} &:& 
\text{"the image of $\mu$ is convex, for any $T^{m}$-action".}
\end{eqnarray*}
Then we have 
\begin{eqnarray*}
(1) &\Leftrightarrow& \text{$A_{m}$ holds for all $m$}, \\
(2) &\Leftrightarrow& \text{$B_{m}$ holds for all $m$}. 
\end{eqnarray*}
At first we see that $A_{1}$ holds by using Proposition \ref{BM} 
and the fact that level sets of a Bott-Morse function on a connected 
compact manifold are 
connected if the critical manifolds all have index and coindex 
$\not= 1$ (see \cite{MS} for example). $B_{1}$ holds clearly 
because in $\mathbb{R}$ connectedness is convexity. 

Now we prove $A_{m-1} \Longrightarrow B_{m}$. Choose a matrix 
$A \in \mathbb{Z}^{m \times (m-1)}$ of maximal rank. 
If we identify $A$ with a linear mapping 
$A: \mathbb{R}^{m-1} \longrightarrow \mathbb{R}^{m}$, then 
$A$ induces an action of $T^{m-1}$-action on $M$ by 
$\theta \cdot p := (A \theta) \cdot p$ for each 
$\theta \in T^{m-1}$ and $p \in M$. This $T^{m-1}$-action 
is a Hamiltonian action with a generalized moment map 
$\mu_{A}(p) := A^{t}\mu(p)$ and a moment one form 
$\alpha_{A}^{\xi} := \alpha^{A\xi}$, where $A^{t}$ denotes the 
transpose of $A$. 

Given any $a \in \mathfrak{t}^{*}$ and $p_{0} \in \mu_{A}^{-1}(a)$, 
\begin{equation*}
p \in \mu_{A}^{-1}(a) \Leftrightarrow 
A^{t}\mu(p) = a = A^{t}\mu(p_{0}) \Leftrightarrow 
\mu(p) - \mu(p_{0}) \in \ker A^{t} 
\end{equation*}
so that 
\begin{equation*}
\mu_{A}^{-1}(a) = \{ p \in M\ |\ \mu(p) - \mu(p_{0}) \in \ker A^{t} \}. 
\end{equation*}
By the statement $A_{m-1}$, $\mu_{A}^{-1}(a)$ is connected. 
Therefore, if we connect $p_{0}$ to $p_{1}$ by a path $p_{t}$ 
in $\mu_{A}^{-1}(a)$, we obtain a path 
$\mu(p_{t}) - \mu(p_{0})$ in $\ker A^{t}$. Since $A^{t}$ is surjective, 
so $\ker A^{t}$ is $1$-dimensional. Hence $\mu(p_{t})$ must go through 
any convex combination of $\mu(p_{0})$ and $\mu(p_{1})$, 
which shows that any point on the line segment from 
$\mu(p_{0})$ to $\mu(p_{1})$ must be in $\mu(M)$. 

Any $p_{0}, p_{1} \in M$ with $\mu(p_{0}) \not= \mu(p_{1})$ can be 
approximated arbitrarily closely by points $p_{0}^{\prime}$ and 
$p_{1}^{\prime}$ with $\mu(p_{1}^{\prime}) - \mu(p_{0}^{\prime}) 
\in \ker A^{t}$ for a matrix $A \in \mathbb{Z}^{m \times (m-1)}$ 
of maximal rank. By the argument above, we see that the line 
segment from $\mu(p_{0}^{\prime})$ to $\mu(p_{1}^{\prime})$ must be in $\mu(M)$. 
By taking limits $p_{0}^{\prime} \longrightarrow p_{0}$, and 
$p_{1}^{\prime} \longrightarrow p_{1}$ we can conclude that $\mu(M)$ is
convex. 

Next we prove $A_{m-1} \Longrightarrow A_{m}$. By identifying 
$\mathfrak{t}$ with $\mathbb{R}^{m}$, we write 
$\mu = (\mu_{1}, \cdots, \mu_{m})$. The generalized moment map $\mu$ 
is called effective if the $1$-forms $d\mu_{1}, \cdots, d\mu_{m}$ 
are linearly independent. Note that $p \in M$ is a regular point 
of $\mu$ if and only if $(d\mu_{1})_{p}, \cdots, (d\mu_{m})_{p}$ 
are linearly independent. 
\begin{lem}
If $\mu$ is not effective, the action reduces to a Hamiltonian 
action of an $(m-1)$-dimensional subtorus. 
\end{lem}
\begin{proof}
If $\mu$ is not effective, there exists 
$0 \not= c = (c_{1}, \cdots, c_{m}) \in \mathbb{R}^{m}$ 
such that $\sum_{i=1}^{m}c_{i}d\mu_{i} = 0$. So if we denote the 
canonical basis of $\mathfrak{t} \cong \mathbb{R}^{m}$ 
by $\xi_{1}, \cdots, \xi_{n}$, then we have 
\begin{equation*}
\sum_{i=1}^{m}c_{i} \Big( (\xi_{i})_{M} + \alpha_{i} \Big)
= \sum_{i=1}^{m}c_{i} 
\Big((\xi_{i})_{M} -\sqrt{-1}( d\mu_{i} + \sqrt{-1}\alpha_{i}) \Big)
\in L, 
\end{equation*}
where $\alpha = (\alpha_{1}, \cdots, \alpha_{m})$. Since 
$\sum_{i=1}^{m}c_{i} \left( (\xi_{i})_{M} + \alpha_{i} \right)$ 
is real and $L \cap \bar L = \{ 0 \}$, we obtain 
$\sum_{i=1}^{m}c_{i}(\xi_{i})_{M} = 0$. 
Consider $\xi = \sum_{i=1}^{m}c_{i}\xi_{i} \in \mathfrak{t}$. 
Then $\mu^{\xi}$ is constant along $M$ because $\xi_{M} = 0$. For the 
simplicity, we may assume $\xi_{1}, \cdots, \xi_{m-1}, \xi$ are 
linearly independent. Then the $T^{m-1}$-action generated by 
$\xi_{1}, \cdots, \xi_{m-1}$ is a Hamiltonian action with 
a generalized moment map $(\mu_{1}, \cdots, \mu_{m-1})$ and 
a moment one form $(\alpha_{1}, \cdots, \alpha_{m-1})$. 
\end{proof}
So we may assume that $\mu$ is effective. 
Then for each $0 \not= \xi \in \mathfrak{t}$, $\mu^{\xi}$ is not 
a constant function. So the critical manifold ${\rm Crit}(\mu^{\xi})$ is 
an even dimensional proper submanifold. Now consider the union of 
critical manifolds $C = \cup_{\eta \not= 0}{\rm Crit}(\mu^{\eta})$. 

\begin{lem}
The union $C$ is a countable union of even dimensional proper 
submanifold, that is, 
\begin{equation*}
C = \cup_{0 \not= \eta \in \mathbb{Z}^{m}}{\rm Crit}(\mu^{\eta}). 
\end{equation*}
\end{lem}
\begin{proof}
We may only show that 
$C \subset \cup_{0 \not= \eta \in \mathbb{Z}^{m}}{\rm Crit}(\mu^{\eta}).$ 
For a rational vector $\xi = \sum_{i=1}^{m}c_{i}\xi_{i}$, 
$\tilde \xi = (\Pi_{i=1}^{m}q_{i})\xi \in \mathbb{Z}^{m}$, where 
$c_{i} = p_{i}/q_{i}$ and $p_{i}, q_{i} \in \mathbb{Z}$. Since 
\begin{equation*}
{\rm Crit}(\mu^{\xi}) = {\rm Fix}(T^{\xi}) = 
{\rm Fix}(T^{\tilde \xi}) = {\rm Crit}(\mu^{\tilde \xi}), 
\end{equation*}
we see that ${\rm Crit}(\mu^{\xi}) \subset \cup_{0 \not= \eta \in \mathbb{Z}^{m}}{\rm Crit}(\mu^{\eta})$. If $\xi$ is irrational, then for rational 
element $\eta \in \mathfrak{t^{\xi}}$ where $\mathfrak{t^{\xi}}$ 
is the Lie algebra of $T^{\xi}$, we have 
${\rm Fix}(T^{\xi}) \subset {\rm Fix}(T^{\eta})$. So we obtain 
\begin{equation*}
{\rm Crit}(\mu^{\xi}) = {\rm Fix}(T^{\xi}) \subset 
{\rm Fix}(T^{\eta}) \subset 
\cup_{0 \not= \eta \in \mathbb{Z}^{m}}{\rm Crit}(\mu^{\eta}). 
\end{equation*}
Hence we have proved 
$C = \cup_{0 \not= \eta \in \mathbb{Z}^{m}}{\rm Crit}(\mu^{\eta})$.
\end{proof}
In particular, $M \setminus C$ is dense subset of $M$. 
In addition, since the condition $p \in M \setminus C$ 
is equivalent to the condition that 
$(d\mu_{1})_{1}, \cdots, (d\mu_{m})_{p}$ are linearly independent, 
we obtain $M \setminus C$ is open dense subset of $M$. 

\begin{lem}
The set of regular values of $\mu$ in $\mu(M)$ is a 
dense subset of $\mu(M)$. 
\end{lem}
\begin{proof}
For each $a = \mu(p) \in \mu(M)$, there exists a sequence 
$\{p_{i}\}_{i = 1}^{\infty} \subset M \setminus C$ 
which satisfies that $\lim_{i \to \infty}p_{i} = p$. Since $p_{i}$ is 
a regular point of $\mu$, $\mu(M)$ contains a neighborhood of 
$\mu(p_{i})$ by implicit function theorem. Moreover there exists 
a regular value $a_{i} \in \mathfrak{t}^{*}$ which is sufficiently 
close to $\mu(p_{i})$ and $\mu^{-1}(a_{i}) \not= \phi$ 
by Sard's theorem. Hence the sequence $\{ a_{i}\}_{i=1}^{\infty}$ 
approximates $a$. 
\end{proof}
By the similar argument, the set of 
$a = (a_{1}, \cdots, a_{m}) \in \mathfrak{t}^{*}$ that 
$(a_{1}, \cdots, a_{m-1})$ is a regular value of 
$(\mu_{1}, \cdots, \mu_{m-1})$ in $\mu(M)$ is also 
a dense subset of $\mu(M)$. Hence, by continuity, to prove that 
$\mu^{-1}(a)$ is connected for every 
$a = (a_{1}, \cdots, a_{m}) \in \mathfrak{t}^{*}$, it suffics to 
prove that $\mu^{-1}(a)$ is connected whenever $(a_{1}, \cdots, a_{m-1})$ 
is a regular value for the reduced generalized moment map 
$(\mu_{1}, \cdots, \mu_{m-1})$. By the induction hypothesis, 
the submanifold $Q = \cap_{i=1}^{m-1}\mu_{i}^{-1}(a_{i})$ is connected 
whenever $(a_{1}, \cdots, a_{m-1})$ is a regular value for 
$(\mu_{1}, \cdots, \mu_{m-1})$. To finish the proof, we need following 
lemma. 
\begin{lem}\label{BM2}
If $(a_{1}, \cdots, a_{m-1})$ is a regular value for 
$(\mu_{1}, \cdots, \mu_{m-1})$, the function 
$\mu_{m} : Q \longrightarrow \mathbb{R}$ is a Bott-Morse function 
of even index and coindex. 
\end{lem}
\begin{proof}
By the hypothesis, $Q$ is a $2n-(m-1)$ dimensional connected submanifold 
of $M$. For each $p \in Q$, the subspace $W$ of the cotangent space 
$T_{p}^{*}M$ generated by $(d\mu_{1})_{p}, \cdots, (d\mu_{m-1})_{p}$ is 
$(m-1)$ dimensional because $p$ is regular. So the tangent space $T_{p}Q$ 
of $Q$ coincises with the annihilator of $W$; 
\begin{equation*}
T_{p}Q = \{ X \in T_{p}M\ |\ f(X) = 0\ (\forall f \in W) \}. 
\end{equation*}
Hence $p \in Q$ is a critical point of $\mu_{m}: Q \longrightarrow \mathbb{R}$ 
if and only if there exists real numbers $c_{1}, \cdots, c_{m-1}$ 
such that 
\begin{equation*}
\sum_{i=1}^{m-1}c_{i}(d\mu_{i})_{p} + (d\mu_{m})_{p} = 0. 
\end{equation*}
This means that $p$ is a critical point of the function 
$\mu^{\xi} : M \longrightarrow \mathbb{R}$, where 
$\xi = (c_{1}, \cdots, c_{m-1}, 1) \in \mathfrak{t} \cong \mathbb{R}^{m}$. 
By Proposition \ref{BM}, $\mu^{\xi}$ is a Bott-Morse function with 
even index and coindex. We shall prove the critical manifold 
${\rm Crit}(\mu^{\xi})$ intersects $Q$ transversally at $p$, that is, 
\begin{equation*}
T_{p}M = T_{p}{\rm Crit}(\mu^{\xi}) + T_{p}Q. 
\end{equation*}
This is equivalent to $W \cap T_{p}^{*}{\rm Crit}(\mu^{\xi}) = \{ 0 \}$. 
To see this, we need to prove $\sum_{i=1}^{m-1}c_{i}(d\mu_{i})_{p}(X) = 0$ 
for any $X \in T_{p}{\rm Crit}(\mu^{\xi})$ and 
$\sum_{i=1}^{m-1}c_{i}(d\mu_{i})_{p} \in W \cap T_{p}^{*}{\rm Crit}(\mu^{\xi})$. This means that the linear functionals 
$(d\mu_{1})_{p}, \cdots, (d\mu_{m-1})_{p}$ remain linearly 
independent when restricted to the subspace $T_{p}{\rm Crit}(\mu^{\xi})$. 
Consider the vector fields 
$\xi_{1}^{+}, \cdots, \xi_{m-1}^{+}$ on $M$ defined by 
\begin{equation*}
d\mu_{i} = \omega_{+}(\xi_{i}^{+}),\quad i=1, \cdots, m-1. 
\end{equation*}
Then $(\xi_{1}^{+})_{p}, \cdots, (\xi_{m-1}^{+})_{p}$ are 
linearly independent on $T_{p}M$ because $p$ is regular. 
So they are also linearly independent on $T_{p}^{*}{\rm Crit}(\mu^{\xi})$. 
Since the $2$-form $\omega_{+}$ is still non-degenerate when 
it is restricted to ${\rm Crit}(\mu^{\xi})$, 
so $(d\mu_{1})_{p}, \cdots, (d\mu_{m-1})_{p}$ are linearly 
independent on $T_{p}^{*}{\rm Crit}(\mu^{\xi})$ and hense 
${\rm Crit}(\mu^{\xi})$ is transverse to $Q$ as claimed. 

This implies that the orthogonal complement of the subspace 
$T_{p}{\rm Crit}(\mu^{\xi})$ is contained in $T_{p}Q$. Hence 
the Hessian of $\mu^{\xi}$ at $p$ is non-degenerate on this space 
with even index and coindex. In other words, 
${\rm Crit}(\mu^{\xi}) \cap Q$ is the critical manifold of 
$\mu^{\xi}|_{Q}$ of even index and coindex. The same holds for 
$\mu_{m}|_{Q}$ since it only differs from $\mu^{\xi}$ by the 
constant $\sum_{i=1}^{m-1}c_{i}a_{i}$. Thus we have proved 
that the function $\mu_{m} : Q \longrightarrow \mathbb{R}$ is 
a Bott-Morse function with even index and coindex. 
\end{proof}
By applying Lemma \ref{BM2}, if $(a_{1}, \cdots, a_{m-1})$ is a 
regular value for $(\mu_{1}, \cdots, \mu_{m-1})$, then the 
level set $\mu_{m}^{-1}(a_{m}) \cap Q = \mu^{-1}(a)$ is connected. 
This shows that $A_{m-1} \Longrightarrow A_{m}$. 

Finally, we shall prove the third claim, that is, the image of the 
generalized moment map $\mu$ is the convex hull of the images 
of the fixed points of the action. 
By Lemma \ref{Fix}, the fixed point set ${\rm Fix}(T^{m})$ of the 
action decomposes into finitely many even dimensional connected 
submanifolds $C_{1}, \cdots, C_{N}$ of $M$. The generalized moment 
map $\mu$ is constant on each of these sets because 
$C_{i} \subset {\rm Crit}(\mu^{\xi})$ for $i = 1, \cdots, N$ and 
any $\xi \in \mathfrak{t}$. Hence there exists 
$a_{1}, \cdots, a_{N} \in \mathfrak{t}^{*}$ such that 
\begin{equation*}
\mu(C_{i}) = \{ a_{i} \},\quad i=1, \cdots, N. 
\end{equation*}
By what we have proved so far the convex hull of the 
points $a_{1}, \cdots, a_{N}$ is contained in $\mu(M)$. 
Conversely, let $a \in \mathfrak{t}^{*}$ be a point which 
is not in the convex hull of $a_{1}, \cdots, a_{N}$. Then 
there exists a vector $\xi \in \mathfrak{t}$ with rationally 
independent components such that 
\begin{equation*}
a_{i}(\xi) < a(\xi),\quad i=1, \cdots, N. 
\end{equation*}
Since the components of $\xi$ are rationally independent, 
we have ${\rm Crit}(\mu^{\xi}) = {\rm Fix}(T^{m})$. Hence 
the function $\mu^{\xi} : M \longrightarrow \mathbb{R}$ 
attains its maximum on one of the sets $C_{1}, \cdots, C_{N}$. 
This implies 
\begin{equation*}
\sup_{p \in M}\mu^{\xi}(p) < a(\xi), 
\end{equation*}
and hence $a \not \in \mu(M)$. This shows that $\mu(M)$ 
is the convex hull of the points $a_{1}, \cdots, a_{N}$ and 
{\rm Theorem A} is proved. \qed

\begin{rem}
If $M$ is an orbifold and $H$ is a closed $3$-form on $M$, 
we can define notions of $H$-twisted generalized complex structures 
of $M$ and Hamiltonian actions of compact Lie groups on an $H$-twisted generalized complex orbifold in usual way. By applying same arguments of 
our proof and Theorem 5.1 in \cite{Ler3}, Theorem A still holds 
when $M$ is a connected compact $H$-twisted generalized complex 
orbifold. Of course, $C_{1}, \cdots, C_{N}$ are connected suborbifolds 
in this case. 

\end{rem}

\begin{rem}
When the manifold is noncompact and the Lie group is non abelian, 
the convexity and connected properties still hold in the sense 
of Theorem 1.1 and Theorem 4.3 in \cite{Ler2}. In view of the works 
of Lerman, Meinrenken, Tolman, and Woodward in \cite{Ler2}, we can 
import notions of symplectic cuts and Cross-section theorem 
in symplectic geometry to generalized complex geometry. 
These techniques tell us that the convexity and connected properties 
still hold in general cases. Detailed proof can be seen in our work 
\cite{YN}. 
\end{rem}


\begin{thebibliography}{99}
\bibitem{At} 
M. F. Atiyah, 
Convexity and commuting hamiltonians, 
Bull. London Math. Soc. 14 (1982), 1--15. 
\bibitem{GHR}
Jr. S. Gates, C. Hull, and M. Rocek, 
Twisted multiplets and new supersymmetric nonlinear
$\sigma$-model, 
Nuclear Phys. B. 248(1) (1984), 157--186.
\bibitem{GS} 
V. Guillemin, and S. Sternberg, 
Convexity properties of the moment mapping I and I\hspace{-.1em}I, 
Invent. Math. 67 (1982), 491--513; Invent. Math. 77 (1984), 533--546. 
\bibitem{Gua} 
M. Gualtieri, 
Generalized complex geometry, 
PhD thesis, Oxford University, 2004, math.DG/0401221 
\bibitem{Hi}
N.Hitchin, 
Generalized Calabi-Yau manifolds, 
Quart. J. Math. 54, 281-308 (2003) 
\bibitem{KT}
A. Kapustin, and A. Tomasiello, 
The general $(2, 2)$ gauged sigma model with three-form flux, 
preprint, hep-th/0610210.
\bibitem{Ler1}
E. Lerman, 
Symplectic cuts, 
Math. Res. Lett. 2 (1995), no. 3, 247--258. 
\bibitem{Ler2}
E. Lerman, E. Meinrenken, S. Tolman, and C. Woodward, 
Nonabelian convexity by symplectic cuts, 
Topology 37 (1998), no. 2, 245--259. 
\bibitem{Ler3} E. Lerman, S. Tolman, 
Hamiltonian torus actions on symplectic orbifolds and toric varieties. 
Trans. Amer. Math. Soc. 349 (1997), no. 10, 4201--4230. 
\bibitem{Lin}
Y. Lin and S. Tolman, 
Symmetries in generalized K{$\rm \ddot a$}hler geometry, 
Comm. in Math. Phys. 268 (2006) 199-222. 
\bibitem{Lin2}
Y. Lin, 
Generalized geometry, equivariant $\overline{\partial}\partial$-lemma, and torus actions, 
the Journal of Geometry and Physics, 57 (2007) 1842-1860. 
\bibitem{MS}
D. McDuff, and D. Salamon, 
Introduction to Symplectic Topology, 
Oxford Mathematical Monographs, Oxford University Press, New York, 1995. 
\bibitem{YN}
Y. Nitta, 
Convexity properties for generalized moment maps I\hspace{-.1em}I, 
preprint. 
\end{thebibliography}
\end{document}